\documentclass[11pt, reqno]{amsart}
\usepackage{ amsthm, amscd, amsfonts, amssymb, graphicx,tikz, color}
\usepackage[all]{xy}
\usepackage{hyperref}

\newcommand{\Frac}{\mbox{\rm Frac\,}}
\newcommand{\C}{{\mathbb C}}
\begin{document}

\title{\bf On enveloping skew fields of some Lie superalgebras}

\author{Jacques Alev}
\address[J.~Alev]{Universit\'e de Reims, Laboratoire de Math\'ematiques,
Moulin de la Housse, B.P. 1039, 51687 Reims cedex 2 (France)}
\email{jacques.alev@univ-reims.fr}
\author{Fran\c cois Dumas}
\address[F.~Dumas]{Universit\'e Blaise Pascal (Clermont-Ferrand 2), Laboratoire de
Math\'ematiques (UMR 6620 - CNRS), B.P. 80026, 63171 Aubi\`ere cedex (France)}
\email{Francois.Dumas@math.univ-bpclermont.fr}

\begin{abstract} We determine the skew fields of fractions of the enveloping algebra
of the Lie superalgebra $\mathfrak{osp}(1,2)$ and of some significant subsuperalgebras 
of the Lie superalgebra $\mathfrak{osp}(1,4)$. We compare the kinds of skew fields
arising from this ``super" context with the Weyl skew fields in the classical
Gelfand-Kirillov property.
\end{abstract}

\date\today

\maketitle

\section*{Introduction}This paper deals with the question of a possible analogue of the Gelfand-Kirillov
property for the enveloping algebras of Lie superalgebras. Let us recall that a finite dimensional complex
Lie algebra $\mathfrak g$ satisfies the Gelfand-Kirillov property when its enveloping skew field,  that
is the skew field of fractions of the enveloping algebra $\mathcal U(\mathfrak g)$, is isomorphic to
a Weyl skew field over a purely transcendental extension of $\C$. A rich literature has developed on this topic
from the seminal work \cite{GK} and we refer to the  papers \cite{AOVdB}, \cite{Premet}
and their bibliographies for an overview on it. 

A natural starting point for the same problem for
a finite dimensional complex Lie superalgebra $\mathfrak g$ is the classification of
the classical simple Lie superalgebras (see \cite{Kac}, \cite{Mus2}) and more precisely the study
of the orthosymplectic Lie superalgebra $\mathfrak{osp}(1,2n)$ since this is the only case in the
classification whose enveloping algebra is a domain (see \cite{AL}, \cite{GorLan}). 
This topic is introduced and discussed by Musson in \cite{Mus1} who proves in particular
that $\Frac(\mathcal U(\mathfrak{osp}(1,2n)))$ is not isomorphic to  
a Weyl skew field over a purely transcendental extension of $\C$ when $n=1$.
We show here that the same is true for any $n$ describing explicitely some classes
of skew fields arising from this context.
The even part $\mathfrak g_{\overline 0}$
of $\mathfrak g=\mathfrak{osp}(1,2n)$ is the Lie algebra $\mathfrak{sp}(2n)$ of the symplectic group,
for which the Gelfand-Kirillov property remains an open question (see \cite{Premet}).
Therefore we concentrate in this exploratory paper on the case of $\mathfrak{osp}(1,2)$
and on some significant subsuperalgebras of $\mathfrak{osp}(1,4)$. We consider in
$\mathfrak{osp}(1,4)$ the Lie subsuperalgebras $\mathfrak n^+,\mathfrak b^+$ and $\mathfrak p^+$ which have as
even parts respectively the nilpotent positive part, the associated Borel subalgebra and 
the associated parabolic subalgebra in the triangular decomposition of the even part $\mathfrak g_{\overline 0}=\mathfrak{sp}(4)$.
Determining their enveloping skew fields is the content of sections 2 and 3 of the paper.

The skew fields appearing in this ``super" context are skew fields of rational functions mixing
classical Weyl relations $xy-yx=1$ and ``fermionic" relations $xy+yx=1$ (or equivalently
$xy=-yx$ up to rational equivalence) between the generators. A noteworthy fact is that
these relations are braided and not necessarily pairwise separable up to isomorphism as in the case
of the classical Weyl skew fields. The main properties of these skew fields,
which already appeared in \cite{AleDum} and \cite{Rich}, are given in section 1.

We end this introduction by a short reminder on the para-Bose definition of the Lie superalgebra $\mathfrak{osp}(1,2n)$ and its
enveloping algebra (see \cite{GanPal}, \cite{PalVdJ}). The basefield is $\C$. We fix an integer $n\geq 1$. 
We have $\mathfrak{osp}(1,2n)=\mathfrak{g}_{\overline 0}\oplus\mathfrak{g}_{\overline 1}$ where
the even part $\mathfrak{g}_{\overline 0}$ is the Lie algebra $\mathfrak{sp}(2n)$ of the symplectic group
and $\mathfrak{g}_{\overline 1}$ is a vector space of dimension $2n$. 
As a Lie superalgebra, $\mathfrak{osp}(1,2n)$ is generated by the  $2n$ elements 
$b_i^{\pm}$ ($1\leq i\leq n$) of a basis of the odd part $\mathfrak{g}_{\overline 1}$.
The $2n^2+n$ elements $\{b_j^\pm,b_k^\pm\}$ ($1\leq j\leq k\leq n$) and $\{b_j^+,b_k^-\}$ ($1\leq j,k\leq n$)
form a basis of $\mathfrak{g}_{\overline 0}$. The dimension of the vector space $\mathfrak{osp}(1,2n)$
is $2n^2+3n$. The brackets are given by the so called ``parabose" relations:
\begin{equation}\label{parabose}
[\{ b_j^\xi, b_k^\eta\}, b_\ell^\epsilon]= (\epsilon-\xi )\delta_{j\ell}b_k^\eta+(\epsilon-\eta)\delta_{k\ell}b_j^\xi
\end{equation}\begin{multline}\label{biparabose}
[\{ b_i^\xi, b_j^\eta\},\{ b_k^\epsilon, b_\ell^\varphi\} ]= 
(\epsilon  - \eta)\delta_{jk}\{b_i^\xi , b_\ell^\varphi\} +
(\epsilon  - \xi)\delta_{ik}\{b_j^\eta , b_\ell^\varphi\} \\ + 
(\varphi  - \eta)\delta_{j\ell}\{b_i^\xi , b_k^\epsilon\} + 
(\varphi  - \xi)\delta_{i\ell}\{b_j^\eta , b_k^\epsilon\}.\hspace{1cm} 
\end{multline}
By the PBW theorem (see \cite{Mus2}), the enveloping algebra 
$\mathcal U(\mathfrak{osp}(1,2n))$ is generated by the $2n^2+n$ elements:
\begin{align}
&b_i^\pm, \  k_i:=\textstyle\frac12\{b_i^-,b_i^+\} \ \ \text{for} \ 1\leq i\leq n,\label{U1}\\
&a_{ij}^\pm:=\textstyle\frac12\{b_{i}^\pm,b_{j}^\pm\},  \
s_{ij}:=\frac12\{b_{i}^-,b_{j}^+\}, \ t_{ij}:=\frac12\{b_{i}^+,b_{j}^-\} \ 
\text{for} \ 1\leq i<j\leq n,\label{U2}
\end{align}
with commutation relations deduced from  \eqref{parabose} and \eqref{biparabose}
taking $\{x,y\}=xy+yx$ if $x,y\in\mathfrak g_{\overline 1}$, and $[x,y]=xy-yx$ otherwise. The enveloping algebra
$\mathcal U(\mathfrak{sp}(2n))$ of the even part is the subalgebra of 
$\mathcal U(\mathfrak{osp}(1,2n))$ generated by the $2n^2+n$ elements 
$(b_i^\pm)^2, k_i$ for $1\leq i\leq n$, and 
$a_{ij}^\pm, s_{ij}, t_{ij}$ for $1\leq i<j\leq n$.


\section{Some skew fields}\label{firstsection}
\subsection{Definitions and notations} We fix the basefield to be $\C$.
As usual $\mathbf{A}_1$ is the Weyl algebra, that is the algebra generated
over $\C$ by two generators $x,y$ satisfying the commutation law $xy-yx=1$.
We also define $\mathbf{A}^1$ as the algebra generated 
over $\C$ by two generators $u,v$ satisfying the commutation law $uv+vu=1$.
For any nonnegative integers $r,s$, we denote by $\mathbf{A}_r^s$ the $\C$-algebra:
\[\mathbf{A}_r^s
=\underbrace{\mathbf{A}_1\otimes\mathbf{A}_1\otimes\cdots\otimes\mathbf{A}_1}_{r \text{ factors}}\otimes
\underbrace{\mathbf{A}^1\otimes\mathbf{A}^1\otimes\cdots\otimes\mathbf{A}^1}_{s \text{ factors}}\]
This is clearly a noetherian domain, so we can define $\mathbf{D}_r^s=\Frac\mathbf{A}_r^s$.
For any integer $t\geq 0$, we consider the polynomial algebra $\C[z_1,\ldots,z_t]$, the noncommutative
noetherian domain $\mathbf A_{r,t}^s=\mathbf A_r^s\otimes\C[z_1,\ldots,z_t]$ and its skew field
of fractions $\mathbf D_{r,t}^s=\Frac\mathbf A_{r,t}^s$. In particular $\mathbf D_{r}^s=\mathbf D_{r,0}^s$.
If $s=0$ and \mbox{$r\geq 1$}, $\mathbf D_{r,t}^0$ is the classical Weyl skew field 
usually denoted by $D_r(\C(z_1,\ldots,z_t))$ or $\mathcal D_{r,t}(\C)$.
The following lemma asserts that the $\mathbf{D}_{r,t}^s$ are a particular case of the mixed Weyl skew fields
considered in \cite{AleDum}.
\subsection{Lemma}\label{lemmaD} {\sl The skew field  $\mathbf{D}_{r,t}^s$ is isomorphic
to the skew field of fractions of the algebra $\widehat{\mathbf{A}}_{r,t}^s$ generated over $\C$ by $2r+2s+t$ 
generators $x_1,\ldots,x_r$, $y_1,\ldots y_r$, $u_1,\ldots,u_s$, $w_1,\ldots,w_s$, $z_1,\ldots,z_t$ 
satisfying the commutation relations:
\begin{align*}
[x_i,y_i]=1, \ [x_i,y_j]=[x_i,x_j]=[y_i,y_j]=0  \ &\quad(1\leq i\not=j\leq r),\\
u_iw_i=-w_iu_i, \ [u_i,w_j]=[u_i,u_j]=[w_i,w_j]=0  \ &\quad(1\leq i\not=j\leq s),\\
[x_i,u_j]=[x_i,w_j]=[y_i,u_j]=[y_i,w_j]=0  \ &\quad(1\leq i\leq r, 1\leq j\leq s)\\
[x_i,z_k]=[y_i,z_k]=[u_j,z_k]=[w_j,z_k]=[z_k,z_\ell]=0  \ &\quad(1\leq i\leq r, 1\leq j\leq s, \\
\  &\quad 1\leq k,\ell\leq t).\end{align*}}\begin{proof}For any 
$1\leq i\leq s$, let us consider the copy of $\mathbf{A}^1$ generated by $u_i,v_i$ with relation
$u_iv_i+v_iu_i=1$. The element $w_i:=u_iv_i-v_iu_i=2u_iv_i-1$ of $\mathbf{A}^1$ satisfies
$w_iu_i=-u_iw_i$ and $w_iv_i=-v_iw_i$. In the skew field of fractions, the subfield generated by $u_i,v_i$
is isomorphic to the subfield generated by $u_i,w_i$ since $v_i=\frac12u_i^{-1}(w_i+1)$. Hence the proof is complete.\end{proof}

We sum up in the following proposition some basic facts about the skew fields $\mathbf{D}^s_{r,t}$.
It shows in particular that for $s\not=0$ the skew fields $\mathbf{D}_{r,t}^s$ are not isomorphic to
classical Weyl skew fields.
\subsection{Proposition}\label{propD}{\sl Let $r,s,t$ be any nonnegative integers. Then:\begin{itemize}
\item[{\rm(i)}] the Gelfand-Kirillov transcendence degree of $\mathbf{D}^s_{r,t}$ equals to $2r+2s+t$;
\item[{\rm(ii)}] the center of $\mathbf{D}^s_{r,t}$ is 
$\C(u_1^2,\ldots,u_s^2,w_1^2,\ldots,w_s^2,z_1,\ldots,z_t)$,  with the notations of lemma \ref{lemmaD};
\item[{\rm(iii)}] $\mathbf{D}^s_{r,t}$ is isomorphic to a classical Weyl skew field
$\mathbf{D}^0_{r',t'}$ if and only if $s=0, r=r'$ and $t=t'$.\end{itemize}}
\begin{proof} The algebra $\widehat{\mathbf{A}}_{r,t}^s$ of lemma \ref{lemmaD} is a particular case
of the algebras $S_{n,r}^\Lambda$ studied in \cite{Rich}. Explicitly $\widehat{\mathbf{A}}_r^s=S_{n,r}^\Lambda$
for $n=r+2s+t$ and $\Lambda=(\lambda_{ij})$ the $n\times n$ matrix with entries in $\C$ defined by 
$\lambda_{r+2k-1,r+2k}=\lambda_{r+2k,r+2k-1}=-1$ for any $1\leq k\leq s$, and $\lambda_{i,j}=1$ in any other case.
Then points (i) and (ii) follow respectively from proposition 1.1.4 and proposition 3.3.1 of \cite{Rich}. 
Suppose now that $\mathbf{D}^s_{r,t}$ is isomorphic to $\mathbf{D}^0_{r',t'}$ for some $r'\geq 1,t'\geq 0$. Denote
$G(L)=(L^\times)'\cap \C^\times$ the trace on $\C^\times$ of the commutator
subgroup of the group of nonzero elements of $L$ for any skew field $L$ over $\C$.
It follows from theorem 3.10 of \cite{AleDum1} that $G(\mathbf{D}^0_{r',t'} )=\{1\}$ while
it is clear by lemma \ref{lemmaD} that $-1\in G(\mathbf{D}^s_{r,t})$ if $s\geq 1$. Hence $s=0$. Then comparing the centers
we deduce $t=t'$ and comparing the Gelfand-Kirillov transcendence degrees we conclude $r=r'$.\end{proof}
\subsection{Remark}Each copy in $\mathbf{D}_{r,t}^s$ of the algebra $\widehat{\mathbf A}^1$ generated
over $\C$ by two generators $u,w$ satisfying $uw=-wu$ can be viewed as the enveloping algebra
of the nilpotent Lie superalgebra $\mathfrak f=\mathfrak f_{\overline 0}\oplus\mathfrak f_{\overline 1}$ where $\mathfrak{f}_{\overline 0}=\C z\oplus\C t$
and $\mathfrak{f}_{\overline 1}=\C u\oplus\C w$ with brackets $\{u,u\}=z$, $\{w,w\}=t$, $\{u,w\}=0$.\medskip

The argument used in the proof of point (iii) of the previous proposition allows
to show the following proposition, as predicted in \cite{Mus1}.

\subsection{Proposition}\label{Fracosp}{\sl For any integer $n\geq 1$,
$\Frac(\mathcal U(\mathfrak{osp}(1,2n)))$ is not isomorphic 
to a classical Weyl skew field $\mathbf{D}^0_{r,t}$ for any $r\geq 1,t\geq 0$. 

More generally, for any subsuperalgebra $\mathfrak g$
of $\mathfrak{osp}(1,2n)$ containing the generators $b_i^+$ and $b_i^-$ for some $1\leq i\leq n$, 
$\Frac(\mathcal U(\mathfrak g))$ is not isomorphic to a classical Weyl skew field $\mathbf{D}^0_{r,t}$
for any $r\geq 1,t\geq 0$.}
\begin{proof}
If $\mathfrak g$ contains $b_i^+$ and $b_i^-$, it contains the element $k_i=\frac12\{b_i^+,b_i^-\}$.
Then $\mathcal U(\mathfrak g)$ contains the element $z_i=b_i^+b_i^--b_i^-b_i^++1=2b_i^+b_i^--2k_i+1$.
Using relation \eqref{parabose}, we have $[k_i,b_i^+]=b_i^+$. An obvious calculation
gives $z_ib_i^+=-b_i^+z_i$. It follows that $-1\in G(\Frac(\mathcal U(\mathfrak g)))$ ; as at the end 
of the proof of proposition \ref{propD} we conclude that $\Frac(\mathcal U(\mathfrak g))$ cannot
be isomorphic to a classical Weyl skew field.\end{proof}

The skew fields $\mathbf{D}^s_{r,t}$ are the most simple and natural way to mix 
classical Weyl skew fields $\mathcal{D}_{r,t}(\C)$ with ``fermionic" relations $uw=-wu$.
However we will see in the following that they are not sufficient to describe
the rational equivalence of enveloping algebras of Lie superalgebras. 
Some ``braided" versions of mixed skew fields are necessary. 
The low dimensional examples useful for
the following results are introduced in \cite{AleDum}.
Their generalization in any dimension are the subject of a systematic study 
in the article \cite{Rich}.  We recall here their definitions and main properties.

\subsection{Definitions and notations}\label{defF3F4} 
Let $\mathbf S_3$ be the algebra generated over $\C$ by three generators $x,y,z$ satisfying:
\[xy-yx=1, \ \ xz=-zx, \ \ yz=-zy.\]
Crossing two copies of $\mathbf S_3$, we define $\mathbf S_4$ as the algebra generated over $\C$ by 
four generators $x_1,x_2,y_1,y_2$ satisfying:{\renewcommand{\arraystretch}{1.5}\begin{equation*}\begin{matrix}
&x_1y_1-y_1x_1=1, \hfill &x_1y_2=-y_2x_1,\hfill &x_1x_2=-x_2x_1 \hfill \\
&x_2y_2-y_2x_2=1, \hfill &x_2y_1=-y_1x_2, \hfill  &y_1y_2=-y_2y_1.\hfil\end{matrix}\end{equation*}}
The algebras $\mathbf S_3$ and $\mathbf S_4$ 
are obviously noetherian domains. We denote $\mathbf F_3=\Frac\mathbf S_3$ and $\mathbf F_4=\Frac\mathbf S_4$.

The algebra $\mathbf S_4$ is the case $n=2$ of the family of quantum Weyl algebras $A_n^{\overline q,\Lambda}$
introduced in \cite{AleDum1} when all nontrivial entries $\lambda_{ij}$ of $\Lambda$ are equal to $-1$ and all entries
$q_i$ of $\overline q$ are equal to $1$. They have been intensively studied (we refer to \cite{GZ} and to
section 1.3.3 of \cite{Rich2} for a survey and references), are simple of center $\C$ and
have the same Hochschild homology and cohomology as the classical Weyl algebra $A_n(\C)$.
A similar study for $\mathbf S_3$ lies in sections 5 and 7 of \cite{Rich2}.

\subsection{Proposition.}\label{propF3F4} {\sl The following holds for the skew fields $\mathbf F_3$ and $\mathbf F_4$:\begin{itemize}
\item[{\rm(i)}] the Gelfand-Kirillov transcendence degrees of $\mathbf F_3$ and $\mathbf F_4$ are 3 and 4 respectively ;
\item[{\rm(ii)}] the center of $\mathbf F_3$ is  $\C(z^2)$,
and the center of $\mathbf F_4$ is $\C$;
\item[{\rm(iii)}] $\mathbf F_3$ and $\mathbf F_4$ are not isomorphic to $\mathbf{D}^s_{r,t}$,
for any $r,s,t\geq 0$.\end{itemize}}
\begin{proof} These properties are proved under slightly different assumptions in section 3 of \cite{AleDum}.
With the notation of \cite{Rich}, we have $\mathbf S_3=S_{2,1}^\Lambda$ 
and $\mathbf S_4=S_{2,2}^\Lambda$ for $\Lambda=\left(\begin{smallmatrix}1&-1\\-1&1\end{smallmatrix}\right)$.
Then points (i) and (ii) follow respectively from proposition 1.1.4 and proposition 3.3.1 of \cite{Rich}.
Suppose that $\mathbf F_3$ is isomorphic to some $\mathbf D_{r,t}^s$.
Comparing the Gelfand-Kirillov transcendence degree, we necessarily have $(r,s,t)=(0,0,3)$, $(1,0,1)$
or $(0,1,1)$. 
The first case is obviously excluded since $\mathbf F_3$ is not commutative.
The second case is impossible because, with the notation $G(L)$ recalled
in the proof of proposition \ref{propD}, we know that $G(\mathbf D_{1,1}^0)=\{1\}$
by theorem 3.10 of \cite{AleDum1},  while it is clear that 
$-1\in G(\mathbf F_3)$. 
The third case is also impossible because, denoting
$E(L)=[L,L]\cap \C$ the trace on $\C$ of the subspace generated
by the commutation brackets for any skew field $L$ over $\C$,
we have $E(\mathbf D_{0,1}^1)=\{0\}$
by proposition 3.9 of \cite{AleDum1}, and 
$E(\mathbf F_3)=\C$ since $D_1(\C)\subset\mathbf F_3$.
Suppose now that $\mathbf F_4$ is isomorphic to $\mathbf D_{r,t}^s$. 
Since the transcendence degree of the center of $\mathbf D_{r,t}^s$ is at least $t$, it follows
from point (ii) that $s=t=0$. Therefore 
$\mathbf F_4$ would be isomorphic to the usual Weyl skew field $\mathbf D_{2,0}^0$.
One more time this is impossible because $G(\mathbf D^0_{2,0})=\{1\}$ and 
$-1\in G(\mathbf F_4)$.\end{proof}

\subsection{Remarks}\label{F3F4D} 
Let us consider the algebra $\widehat{\mathbf A}_{r,0}^s$
with the notations of lemma \ref{lemmaD}. If $r\geq 1$ and $s\geq 1$, the subfield of $\mathbf D_r^s$ generated by $x_1w_1,y_1w_1^{-1 }$ and $u_1$ is isomorphic to $\mathbf F_3$.
If $r\geq 2$ and $s\geq 1$, the subfield of $\mathbf D_r^s$
generated by $x_1w_1,y_1w_1^{-1 },x_2u_1$ and $y_2u_1^{-1}$ is isomorphic to $\mathbf F_4$.
In other words, $\mathbf F_3$ can be embedded in any $\mathbf{D}^s_r$ such that $r\geq 1,s\geq 1$ and $\mathbf F_4$ can be embedded in  any $\mathbf{D}^s_r$ such that $r\geq 2,s\geq 1$.
More deeply it follows from proposition 5.3.3 of \cite{Rich} that $\mathbf F_4$ cannot
be embedded in some $\mathbf D^s_{r,t}$ for $r\leq 1$.

\subsection{Illustration}\label{graph} We illustrate the definitions of the skew fields under consideration by the following
graphs, stressing the particular nature of the relevant relations. The vertices are parametrized by some system of generators. A directed edge 
$\xymatrix{a \bullet  \ar@*{[|<2pt>]}[r]&  {\bullet}\,b}$ between two generators $a$ and $b$
means that $ab-ba=1$, an undirected edge 
$\xymatrix{a \bullet  \ar@{..}[r]&  {\bullet}\,b}$ means that $ab=-ba$,
and no edge between two generators means that they commute.\bigskip

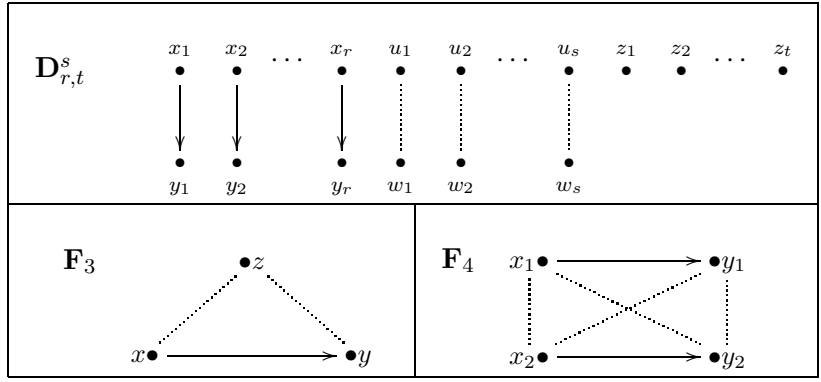
\centerline{\renewcommand{\arraystretch}{2.6}
\begin{tabular}{|c|}\hline
\begin{tabular}{c} $\mathbf D_{r,t}^s$ \\ \end{tabular}\qquad
$\xymatrix{{x_1}\atop{\bullet}\ar@*{[|<2pt>]}[d]\\ {\bullet}\atop{y_1}}$
$\xymatrix{{x_2}\atop{\bullet}\ar@*{[|<2pt>]}[d]\\ {\bullet}\atop{y_2}}$
$\cdots$
$\xymatrix{{x_r}\atop{\bullet}\ar@*{[|<2pt>]}[d]\\ {\bullet}\atop{y_r}}$
$\xymatrix{{u_1}\atop{\bullet}\ar@{..}[d]\\ {\bullet}\atop{w_1}}$
$\xymatrix{{u_2}\atop{\bullet}\ar@{..}[d]\\ {\bullet}\atop{w_2}}$
$\cdots$
$\xymatrix{{u_s}\atop{\bullet}\ar@{..}[d]\\ {\bullet}\atop{w_s}}$
$\xymatrix{{{z_1}\atop{\bullet}}\\}$
$\xymatrix{{{z_2}\atop{\bullet}}\\}$
$\cdots$
$\xymatrix{{{z_t}\atop{\bullet}}\\}$
 \\ \hline
\begin{tabular}{c|c}
\begin{tabular}{cc}
$\mathbf F_3$ &{\small $\xymatrix{
 &   \ar@{..}[ld] \bullet z\ar@{..}[rd]& \\
x \bullet  \ar@*{[|<2pt>]}[rr]& & {\bullet}y}$ }
\end{tabular}
&\begin{tabular}{cc}
$\mathbf F_4$ &{\small $\xymatrix{
x_1 \bullet \ar[rr] & & \bullet y_1\\
x_2 \bullet\ar@{..}[rru] \ar[rr] \ar@{..}[u] & &\ar@{..}[llu]  \ar@{..}[u]  \bullet y_2}$ }
\end{tabular}
\end{tabular}\\ \hline
\end{tabular} }


\section{The enveloping skew field of the Lie superalgebra $\mathfrak{osp}(1,2)$}\label{secondsection}
\subsection{Notations} Applying \eqref{parabose} and \eqref{U1} for $n=1$, the algebra $\mathcal U(\mathfrak{osp}(1,2))$ 
is generated by $b^+, b^-, k$ with relations:\begin{equation}\label{eqosp12}
kb^+-b^+k=b^+,\quad kb^--b^-k=-b^-,\quad b^-b^+=-b^+b^-+2k.\end{equation}
It is clearly an iterated Ore extension 
$\mathcal U(\mathfrak{osp}(1,2))=\C[b^+][k\,;\,\delta][b^-\,;\tau,d]$,
where $\delta$ is the derivation $b^+\partial_{b^+}$ in $\C[b^+]$, $\tau$ is the
automorphism of $\C[b^+][k\,;\,\delta]$ defined by $\tau(b^+)=-b^+$ and $\tau(k)=k+1$, and 
$d$ is the $\tau$-derivation of $\C[b^+][k\,;\,\delta]$ defined by $d(b^+)=2k$ and $d(k)=0$.
\subsection{Proposition}\label{mainresultosp12} {\sl $\Frac\mathcal U(\mathfrak{osp}(1,2))$ is isomorphic
to $\mathbf F_3$.}
\begin{proof} By obvious calculations using \eqref{eqosp12}, the element $z:=b^+b^--b^-b^++1=2b^+b^--2k+1$
satisfies $zb^+=-b^+z$ and $zk=kz$. Since $b^-=\frac12(b^+)^{-1}(z+2k-1)$ in the algebra 
$\mathcal U':=\C(b^+)[k\,;\,\delta][b^-\,;\tau,d]$, we have
$\mathcal U'=\C(b^+)[k\,;\,\delta][z\,;\tau']$ with $kb^+-b^+k=b^+$, $zk=kz$ and $zb^+=-b^+z$. 
Setting  $y:=(b^+)^{-1}k$, we obtain $\mathcal U'=\C(b^+)[y;\,\partial_{b^+}][z\,;\,\tau']$ with
$yb^+-b^+y=1$, $zb^+=-b^+z$ and $zy=-yz$. Hence $\Frac\mathcal U'=\Frac\mathcal U(\mathfrak{osp}(1,2))$
is isomorphic to $\mathbf F_3$.\end{proof}
\subsection{Remarks}
We know that $\mathcal U(\mathfrak{sl}(2))$ is the subalgebra of $\mathcal U(\mathfrak{osp}(1,2))$
generated by $(b^+)^2,(b^-)^2$ and $k$. Actually up to a change of notations  $e:=\frac12(b^+)^2$ and 
$f:=-\frac12(b^-)^2$ it follows from \eqref{eqosp12} that $[k,e]=2e$, $[k,f]=-2f$ et $[e,f]=k$.
We introduce $\omega:=4ef+k^2-2k$ the usual Casimir in $\mathcal U(\mathfrak{sl}(2))$.\medskip

(i) With the notations used in the proof of the previous proposition, we have in $\mathcal U'$ the identities
$f=\frac14e^{-1}(\omega-k^2+2k)$ and $[\frac12e^{-1}k,e]=1$. Therefore $\Frac\mathcal U(\mathfrak{sl}(2))$
is the subfield of $\Frac\mathcal U(\mathfrak{osp}(1,2))$ generated by $e=\frac12(b^+)^2$, $y':=(b^+)^{-2}k=(b^+)^{-1}y$
and $\omega$ with relations $y'e-ey'=1$, $\omega e=e\omega$ and $\omega y'=y'\omega$.
We recover the well known Gelfand-Kirillov property that $\Frac\mathcal U(\mathfrak{sl}(2))$
is a classical Weyl skew field $D_1$ over a center $\C(\omega)$ of transcendence degree one. With the conventions of \ref{graph}, we can illustrate this skew fields embedding by:\medskip

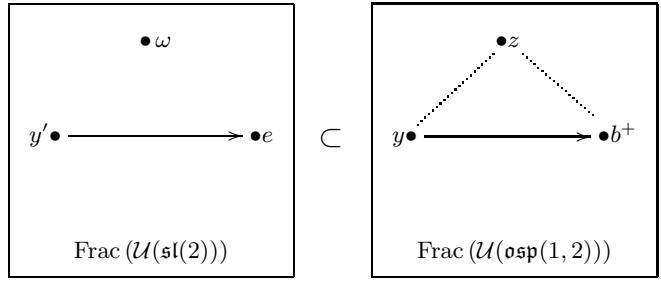
\centerline{\footnotesize{\renewcommand{\arraystretch}{2.1}\begin{tabular}{|c|}\hline
$\xymatrix{  &    \bullet\,\omega & \\
y'\bullet  \ar@*{[|<2pt>]}[rr]& & {\bullet}e}$ \\ \\ $\Frac(\mathcal U(\mathfrak{sl}(2)))$\\
\hline\end{tabular} { \ \large $\subset$ \ }
\begin{tabular}{|c|}\hline
$\xymatrix{  &   \ar@{..}[ld] \bullet z\ar@{..}[rd]& \\
y \bullet  \ar@*{[|<2pt>]}[rr]& & {\bullet}b^+}$ \\ \\
$\Frac(\mathcal U(\mathfrak{osp}(1,2)))$\\ \hline\end{tabular}}}\bigskip

(ii) By previous proposition \ref{mainresultosp12} and point (ii) of proposition \ref{propF3F4},
the center of $\Frac(\mathcal U(\mathfrak{osp}(1,2)))$ is $\C(z^2)$. The element $z$ lying in $\mathcal U(\mathfrak{osp}(1,2))$,
it follows that the center of $\mathcal U(\mathfrak{osp}(1,2))$ is $\C[z^2]$.
A straightforward calculation shows that $z^2=4\omega-2z+3=4\omega-2(z-1)+1$, or equivalently
$(z+1)^2=4(\omega+1)$. Since $z-1=b^+b^--b^-b^+$ by definition of $z$, we recover the well known 
property, see \cite{Pin}, that the center of $\mathcal U(\mathfrak{osp}(1,2))$ is $\C[\theta]$
for $\theta$ the super Casimir operator\begin{equation}\label{Casosp12}
\theta:=\omega -\textstyle\frac12(b^+b^--b^-b^+),\end{equation} 
with $\omega$ the usual Casimir operator of the even part $\mathcal U(\mathfrak{sl}(2))$.
On one hand the above expression of $z^2$ becomes $z^2=4\theta+1$. 
On the other hand, \eqref{Casosp12} implies $z-1=2\omega-2\theta$. 
We deduce that $(2\omega -2\theta+1)^2=4\theta+1$, or equivalently:
\begin{equation}\omega^2-(2\theta-1)\omega+\theta(\theta-2)=0.\end{equation} 
This relation of algebraic dependance between $\theta$ and $\omega$
is exactly the one given in proposition 1.2 of \cite{Pin} up to a normalization 
of the coefficients.\medskip


\section{Enveloping skew fields \\ of some Lie subsuperalgebras of $\mathfrak{osp}(1,4)$}\label{thirdsection}
\subsection{Definitions and notations} We apply for $n=2$ the description of $\mathfrak{osp}(1,2n)$ recalled at the end
of the introduction. We have $\mathfrak{osp}(1,4)=\mathfrak g_{\overline 0}\oplus\mathfrak g_{\overline 1}$ where 
$\mathfrak g_{\overline 1}$ is a vector space of dimension 4 with basis $b_1^+,b_2^+,b_1^-,b_2^-$ and
$\mathfrak g_{\overline 0}$ is the Lie algebra $\mathfrak{sp}(4)$ of dimension 10 with basis:
{\renewcommand{\arraystretch}{1.5}\begin{equation}\label{gen}\begin{matrix}
\textstyle 
c_1^+=\frac12\{b_1^+,b_1^+\},\hfill &c_2^+=\frac12\{b_2^+,b_2^+\},\hfill  
&c_1^-=\frac12\{b_1^-,b_1^-\},\hfill  &c_2^-=\frac12\{b_2^-,b_2^-\},\hfill  \\
a^+=\frac12\{b_1^+,b_2^+\},\hfill  &a^-=\frac12\{b_1^-,b_2^-\},\hfill  
&s=\frac12\{b_1^-,b_2^+\},\hfill  &t=\frac12\{b_1^+,b_2^-\},\hfill  \\
k_1=\frac12\{b_1^-,b_1^+\},\hfill  &k_2=\frac12\{b_2^-,b_2^+\}.\hfill  & &\\
\end{matrix}\end{equation}}
The brackets between these 14 generators of $\mathfrak{osp}(1,4)$ are computed by
the relations \eqref{parabose} et \eqref{biparabose}. 
By \eqref{U1} and \eqref{U2} the algebra $\mathcal U(\mathfrak{osp}(1,4))$ is 
generated by the 10 elements $b_1^+, b_2^+, b_1^-, b_2^-,k_1,k_2,a^+,a^-, s,t$.
The enveloping algebra $\mathcal U(\mathfrak{sp}(4))$ of the even part of $\mathfrak{osp}(1,4)$
is the subalgebra generated by $(b_1^+)^2, (b_2^+)^2,(b_1^-)^2, (b_2^-)^2,k_1,k_2,a^+,a^-, s,t$.\medskip

We describe now some subsuperalgebras  of $\mathfrak{osp}(1,4)$ whose enveloping skew field we study in the following. The even part of each of them
satisfies the usual Gelfand-Kirillov property.
\subsubsection{The nilpotent subsuperalgebra $\mathfrak n^+$}\label{N}  
We define in $\mathfrak g_{\overline 1}$ the subspace $\mathfrak g_{\overline 1}^+:=\C b_1^+\oplus \C b_2^+$
and in $\mathfrak g_{\overline 0}$ the subspace
$\mathfrak n_{\overline 0}^+:=\C c_1^+\oplus \C c_2^+\oplus \C a^+\oplus \C t$. We denote
$\mathfrak n^+:=\mathfrak n_{\overline 0}^+\oplus\mathfrak g_{\overline 1}^+$. We calculate in $\mathfrak{osp}(1,4)$ the 17 brackets 
between the 6 generators of $\mathfrak n^+$ :{\renewcommand{\arraystretch}{1.5}\begin{equation}\label{bracN+}\begin{matrix}
\textstyle 
[a^+,c_1^+]=0,\hfill &[a^+,c_2^+]=0,\hfill &[t,c_1^+]=0,\hfill &[t,c_2^+]=2a^+,\hfill \\
[c_1^+,c_2^+]=0,\hfill & [t,a^+]=c_1^+,\hfill & &\\
\{b_1^+ , b_1^+\}= 2c_1^+,\hfill & \{ b_2^+,b_2^+ \}= 2c_2^+,\hfill &\{b_1^+ ,b_2^+ \}= 2a^+,\hfill& \\
[t,b_1^+]=0,\hfill &[t,b_2^+]=b_1^+,\hfill &[a^+,b_1^+]=0 ,\hfill& [a^+,b_2^+]=0 ,\hfill &\\
[b_1^+,c_1^+]= 0,\hfill &[b_1^+,c_2^+]=0 ,\hfill &[b_2^+,c_1^+]= 0,\hfill &[b_2^+,c_2^+]=0 .\hfill
\end{matrix}\end{equation}}It follows that $\mathfrak n^+$ is a Lie subsuperalgebra of $\mathfrak{osp}(1,4)$ and that
$\mathfrak n_{\overline 0}^+$ is a Lie subalgebra of $\mathfrak g_{\overline 0}$. 
Moreover setting\begin{equation}\label{chev}
x_1:=t,  \ x_2:=c_2^+, \ x_3:=2a^+, \ x_4:=2c_1^+,\end{equation} we rewrite the relations of the first two rows of \eqref{bracN+} as:
{\renewcommand{\arraystretch}{1.5}\begin{equation*}\label{genchevN+}\begin{matrix}
&[x_1 ,x_2 ] = x_3, \hfill &[x_1 ,x_3 ] = x_4,    \hfill   &[x_2 ,x_3 ] = 0,  \hfill   & \\
& [x_1 ,x_4 ] = 0,    \hfill  &[x_2 ,x_4 ] = 0,  \hfill   & [x_3 ,x_4 ] = 0,\hfill         \\
\end{matrix}\end{equation*}}which are the relations between the Chevalley
generators in the enveloping algebra of the nilpotent positive part corresponding
to the root system of type $B_2$. 
We conclude that in the Lie subsuperalgebra  $\mathfrak n^+$ of $\mathfrak g=\mathfrak{osp}(1,4)$, 
the even part $\mathfrak n_{\overline 0}^+$ is isomorphic to the nilpotent positive part in the 
triangular decomposition of $\mathfrak g_{\overline 0}=\mathfrak{sp}(4)$.
\subsubsection{The Borel subsuperalgebra $\mathfrak b^+$}\label{B} We still denote $\mathfrak g_{\overline 1}^+=\C b_1^+\oplus \C b_2^+$
and we introduce in $\mathfrak g_{\overline 0}$ the subspaces 
$\mathfrak h:=\C k_1\oplus\C k_2$ and $\mathfrak b_{\overline 0}^+:=\mathfrak n_{\overline 0}^+\oplus\mathfrak h$.
We define $\mathfrak b^+:=\mathfrak b_{\overline 0}^+\oplus\mathfrak g_{\overline 1}^+$.  
We calculate in $\mathfrak{osp}(1,4)$ the 30 brackets 
between the 8 generators of $\mathfrak b^+$, adding to the 17 brackets of \eqref{bracN+}
the 13 brackets related to $k_1,k_2$, that is :
{\renewcommand{\arraystretch}{1.5}\begin{equation}\label{bracB+}\begin{matrix}
[k_1,b_1^+]= b_1^+,\hfill & \ [k_1 ,b_2^+]= 0,\hfill & \ [k_2,b_1^+]= 0 ,\hfill & \ [k_2,b_2^+]= b_2^+,\hfill \\
[k_1,c_1^+]= 2c_1^+,\hfill & \ [k_1 ,c_2^+]= 0,\hfill & \ [k_2,c_1^+]= 0, \hfill & \ [k_2,c_2^+]= 2c_2^+,\hfill \\
[k_1,a^+]= a^+,\hfill & \ [k_1 ,t]= t,\hfill & \ [k_2,a^+]=a^+ ,\hfill & \ [k_2 ,t]= -t,\hfill \\
[k_1,k_2]= 0.\hfill & & & \\
\end{matrix}\end{equation}}It follows that $\mathfrak b^+$ is a Lie subsuperalgebra of $\mathfrak{osp}(1,4)$ and that
$\mathfrak b_{\overline 0}^+$ is a Lie subalgebra of $\mathfrak g_{\overline 0}$ containing as direct summands the nilpotent Lie subalgebra
$\mathfrak n_{\overline 0}^+$ and the abelian Lie subalgebra $\mathfrak h$. The change of basis \eqref{chev} in $\mathfrak n_{\overline 0}^+$
and the change of basis\begin{equation}\label{cartan}
h_1:=k_2,\qquad h_2:=k_1-k_2,\end{equation} in $\mathfrak h$ allow to rewrite the action of $\mathfrak h$ on $\mathfrak n_{\overline 0}^+$ as:
{\renewcommand{\arraystretch}{1.5}\begin{equation}\label{genchevB+}\begin{matrix}
\textstyle 
[h_1,x_1]= -x_1,\hfill & \ [h_1 ,x_2]= 2x_2,\hfill & \ [h_1,x_3]= x_3, \hfill & \ [h_1,x_4]= 0,\hfill \\
[h_2,x_1]= 2x_1,\hfill & \ [h_2 ,x_2]= -2x_2,\hfill & \ [h_2,x_3]=0 ,\hfill & \ [h_2 ,x_4]= 2x_4.\hfill \\
\end{matrix}\end{equation}} We conclude that 
in the Lie subsuperalgebra $\mathfrak b^+$ of $\mathfrak g=\mathfrak{osp}(1,4)$, the even part 
$\mathfrak b_{\overline 0}^+$ is isomorphic to the positive Borel subalgebra in the triangular decomposition of
$\mathfrak g_{\overline 0}=\mathfrak{sp}(4)$,
and the abelian Lie subalgebra $\mathfrak h$ is isomorphic to the corresponding Cartan subalgebra.
\subsubsection{The parabolic subsuperalgebra $\mathfrak p^+$}\label{P} We introduce in the odd part $\mathfrak g_{\overline 1}$ 
of $\mathfrak{osp}(1,4)$ the subspace $\mathfrak p_{\overline 1}^+:=\mathfrak g_{\overline 1}^+\oplus\C b_2^-=\C b_1^+\oplus \C b_2^+\oplus\C b_2^-$
and in the even part $\mathfrak g_{\overline 0}$ the subspace
$\mathfrak p_{\overline 0}^+:=\mathfrak b_{\overline 0}^+\oplus\C c_2^-=\mathfrak n_{\overline 0}^+\oplus\mathfrak h\oplus\C c_2^-$.
We define $\mathfrak p^+:=\mathfrak p_{\overline 0}^+\oplus\mathfrak p_{\overline 1}^+$.
We calculate in $\mathfrak{osp}(1,4)$ the 48 brackets 
between the 10 generators of $\mathfrak p^+$, adding to the 30 brackets of \eqref{bracN+} and \eqref{bracB+}
the 18 brackets related to  $b_2^-,c_2^-$, that is :
{\renewcommand{\arraystretch}{1.5}\begin{equation}\label{bracP+}\begin{matrix}
\{b_2^-,b_2^+\}=2k_2,\hfill &\qquad [b_2^-, k_1 ] = 0,       \hfill   &\qquad [b_2^-, a^+ ] = b_1^+,\hfill \\
\{b_2^-,b_1^+\}=2t,    \hfill &\qquad [b_2^-, k_2 ] = b_2^-,   \hfill   &\qquad [b_2^-, t ] = 0,\hfill\\ 
[c_2^-,c_2^+]= 4k_2,\hfill &\qquad [c_2^-, k_1 ] = 0,       \hfill   &\qquad [c_2^-, a^+ ] = 2t,\hfill \\
[c_2^-,c_1^+]= 0,    \hfill &\qquad [c_2^-, k_2 ] = 2c_2^-,   \hfill   &\qquad [c_2^-, t ] = 0,\hfill\\ 
[c_2^-,b_1^+]= 0,    \hfill &\qquad [c_2^-, b_2^+ ] = 2b_2^-,   \hfill   &\qquad [c_2^-, b_2^- ] = 0,\hfill\\ 
\{b_2^-,b_2^-\}=c_2^-,\hfill &\qquad [b_2^-, c_1^+] = 0.\hfill & \qquad [b_2^-,c_2^+]=2b_2^+.\hfill\\ \end{matrix}\end{equation}}
It follows that $\mathfrak p^+$ is a Lie subsuperalgebra of $\mathfrak{osp}(1,4)$ and that
$\mathfrak p_{\overline 0}^+$ is a Lie subalgebra of $\mathfrak g_{\overline 0}$ containing as direct summands the Borel subalgebra
$\mathfrak b_{\overline 0}^+$ and the line  $\C c_2^-$. The changes of basis \eqref{chev} and \eqref{cartan} 
allow to rewrite the action of $c_2^-$ on $\mathfrak b_{\overline 0}^+$ as:
{\renewcommand{\arraystretch}{1.5}\begin{equation}\label{genchevP+}\begin{matrix}
\textstyle [c_2^-,x_1]= 4h_1,\hfill & \ [c_2^- ,x_2]= 0 ,\hfill & \ [c_2^-,x_3]= -4x_2, \hfill & \ [c_2^-,x_4]= 0,\hfill \\
[c_2^-,h_1]= 2c_2^-,\hfill & \ [c_2^- ,h_2]= -2c_2^-.\hfill &  & \hfill \\ \end{matrix}\end{equation}}
We conclude that 
in the Lie subsuperalgebra $\mathfrak p^+$ of $\mathfrak g=\mathfrak{osp}(1,4)$, the even part 
$\mathfrak p_{\overline 0}^+$ is isomorphic to the positive parabolic subalgebra in the triangular decomposition of $\mathfrak g_{\overline 0}=\mathfrak{sp}(4)$.
\subsubsection{Remark: the Levi subsuperalgebra $\mathfrak l$ associated to $\mathfrak p^+$}\label{L}
It follows from relations \eqref{bracB+} and \eqref{bracP+} that the subspace $\mathfrak l:=\mathfrak l_{\overline 0}\oplus\mathfrak l_{\overline 1}$
with $\mathfrak l_{\overline 1}:=\C b_2^+\oplus\C b_2^-$ in $\mathfrak p_{\overline 1}^+$ and $\mathfrak l_{\overline 0}:=\C c_2^+\oplus \C k_2\oplus \C c_2^-$
in  $\mathfrak p_{\overline 0}^+$ is a Lie subsuperalgebra of $\mathfrak{osp}(1,4)$. It is clear
that $\mathfrak l$ is isomorphic to $\mathfrak{osp}(1,2)$. The Lie algebra  
$\mathfrak l_{\overline 0}$ is the Levi subalgebra associated to $\mathfrak p_{\overline 0}^+$ in $\mathfrak g_{\overline 0}$
and is isomorphic to $\mathfrak{sl}(2)$.

\subsection{Proposition.}\label{propN} {\sl $\Frac\mathcal U(\mathfrak n^+)$ is isomorphic to 
$\Frac(\mathbf A_1\otimes\mathbf A^1)=\mathbf D_1^1$.}
\begin{proof} By \ref{N}, $\mathcal U(\mathfrak n^+)$ is an
iterated Ore extension $\C[b_1^+,a^+][b_2^+\,;\,\tau,d][t\,;\,\delta]$ expressing the commutation relations
{\renewcommand{\arraystretch}{1.5}\begin{equation}\label{eqVn+}\begin{matrix}
&b_1^+a^+=a^+b_1^+, \hfill &\quad b_2^+a^+=a^+b_2^+,    \hfill   &\quad b_2^+b_1^+=-b_1^+b_2^++2a^+, \hfill \\
&tb_1^+=b_1^+t,    \hfill  &\quad ta^+=a^+t+(b_1^+)^2,  \hfill   &\quad tb_2^+=b_2^+t+b_1^+.         \hfill         \\ 
\end{matrix}\end{equation}}This is a particular case of the more
general theorem 2.1 of \cite{Mus}.
In the algebra $\mathcal U':=\C(b_1^+)[a^+][b_2^+\,;\,\tau,d][t\,;\,\delta]$,
the elements:\begin{equation}\label{t'y}
t':=(b_1^+)^{-2}t,\ \ y:=\frac12(b_1^+b_2^+-b_2^+b_1^+)=b_1^+b_2^+-a^+\end{equation}
satisfy $\mathcal U'=\C(b_1^+)[a^+][y\,;\,\tau][t'\,;\,\delta']$ with relations:
{\renewcommand{\arraystretch}{1.5}\begin{equation}\label{relV'n+}\begin{matrix}
&b_1^+a^+=a^+b_1^+, \hfill &\quad ya^+=a^+y,    \hfill   &\quad yb_1^+=-b_1^+y, \hfill \\
&t'b_1^+=b_1^+t',    \hfill  &\quad t'y=yt',  \hfill   &\quad t'a^+-a^+t'=1.         \hfill         \\ 
\end{matrix}\end{equation}}Hence by lemma \ref{lemmaD} we conclude that
$\Frac\mathcal U(\mathfrak n^+)=\Frac\mathcal U'$ is isomorphic to $\mathbf D_1^1$.\end{proof}

\subsection{Remark} The enveloping algebra $\mathcal U(\mathfrak n_{\overline 0}^+)$ is the subalgebra of $\mathcal U(\mathfrak n^+)$ 
generated by $(b_1^+)^2,(b_2^+)^2,a^+, t$ with commutation relations coming from \eqref{bracN+}:
{\renewcommand{\arraystretch}{1.5}\begin{equation*}\label{eqUn+}\begin{matrix}
&(b_1^+)^2a^+=a^+(b_1^+)^2, \hfill &\  (b_2^+)^2a^+=a^+(b_2^+)^2,  \hfill   &\ (b_2^+)^2(b_1^+)^2=(b_1^+)^2(b_2^+)^2,\hfill \\
&t(b_1^+)^2=(b_1^+)^2t,    \hfill  &\  ta^+=a^+t+(b_1^+)^2,  \hfill         &\ t(b_2^+)^2=(b_2^+)^2t+2a^+. \hfill         \\ 
\end{matrix}\end{equation*}}The element $t'$ defined in \eqref{t'y} lies in $\Frac\mathcal U(\mathfrak n_{\overline 0}^+)$
and the element $y$ defined in \eqref{t'y} satisfies $y^2=(a^+)^2-(b_1^+)^2(b_2^+)^2$ which also lies in  $\Frac\mathcal U(\mathfrak n_{\overline 0}^+)$.
Hence $\Frac\mathcal U(\mathfrak n_{\overline 0}^+)$ is the subfield of $\Frac\mathcal U(\mathfrak n^+)$ generated
by $(b_1^+)^2,y,a^+, t'$ with more simple commutation relations:{\renewcommand{\arraystretch}{1.5}\begin{equation*}\label{eqUn+bis}\begin{matrix}
&(b_1^+)^2a^+=a^+(b_1^+)^2, \hfill &\  y^2a^+=a^+y^2,  \hfill   &\ y^2(b_1^+)^2=(b_1^+)^2y^2,\hfill \\
&t'(b_1^+)^2=(b_1^+)^2t',    \hfill  &\  t'a^+-a^+t'=1,  \hfill         &\ t'y^2=y^2t'. \hfill         \\ 
\end{matrix}\end{equation*}}We recover the well known Gelfand-Kirillov property that $\Frac\mathcal U(\mathfrak n_{\overline 0}^+)$
is a classical Weyl skew field $D_1$ over a center $\C((b_1^+)^2,y^2)$ of transcendence degree two.
We have $(b_1^+)^2=\frac12x_4$ et $y^2=\frac14(x_3^2-2x_2x_4$) with notations \eqref{chev}.
Up to a normalization we recover the well known expressions for the generators of
the center of  $\mathcal U(\mathfrak n_{\overline 0}^+)$ in terms of Chevalley generators.\bigskip

The following theorem gives a decomposition of $\Frac\mathcal U(\mathfrak b^+)$
into two commuting subfields respectively isomorphic to  $\mathbf D_1^0$ and $\mathbf F_4$.\medskip

\subsection{Theorem}\label{thmB} {\sl $\Frac\mathcal U(\mathfrak b^+)$ is isomorphic to 
$\Frac(\mathbf A_1\otimes\mathbf S_4)$.}

\begin{proof} By \ref{B}, $\mathcal U(\mathfrak b^+)$ is generated in $\mathcal U(\mathfrak{osp}(1,4))$
by $\mathcal U(\mathfrak n^+)$ and $\mathcal U(\mathfrak h)$ with the commutation relations \eqref{eqVn+} 
and the action of $k_1$ and $k_2$ on $b_1^+, b_2^+, a^+, t$ coming from \eqref{bracB+}. 
Taking again the notations used in the proof of proposition \ref{propN}, this action extends to 
$\mathcal U''=\C(b_1^+,a^+)[y\,;\,\tau][t'\,;\,\delta']$ by:
{\renewcommand{\arraystretch}{1.5}\begin{equation*}\label{kbyat'}\begin{matrix}
&[k_1,b_1^+ ] = b_1^+,\hfill &\quad [k_1,y ] = y,       \hfill   &\quad [k_1, a^+] = a^+,\hfill &\quad [k_1,t' ] =- t',\hfill\\
&[k_2,b_1^+ ] = 0,    \hfill &\quad [k_2,y ] = y,   \hfill   &\quad [k_2, a^+] = a^+,\hfill &\quad [k_2,t'] = -t'.\hfill\\ 
\end{matrix}\end{equation*}}The change of variables:\begin{equation}\label{defk'}k_1'=(b_1^+)^{-1}(k_1-k_2),
\qquad k'_2=(a^+)^{-1}k_2
\end{equation}gives:\begin{equation}\label{k'ba}
k'_1k'_2=k'_2k'_1, \ [k'_1,b_1^+ ] = 1, \ [k'_1, a^+] = 0, \ 
[k'_2,b_1^+ ] = 0,  \ [k'_2, a^+] = 1.\end{equation}That shows that the subalgebra $\mathcal W$  generated
by $b_1^+,a^+,k'_1,k'_2$ in 
$\Frac\mathcal U(\mathfrak b^+)$ is isomorphic to the Weyl algebra $A_2=\mathbf A_1\otimes\mathbf A_1$. The commutation relations
of these new generators $k'_1,k'_2$ with the generator $y$ 
are $k'_1y=-yk'_1$ and $k'_2y=yk'_2+y(a^+)^{-1}$. We replace $y$ by:
\begin{equation}\label{defy'}y':=(a^+)^{-1}y=(a^+)^{-1}b_1^+{b_2^+}-1,\end{equation}which satisfies:
\begin{equation}\label{bracy'}y'k'_2 =k'_2y', \ \ y'a^+=a^+y',\ \ y'b_1^+=-b_1^+y',\ \ y'k'_1=-k'_1y'.\end{equation}
We deduce with \eqref{k'ba} that the subalgebra $\mathcal V$ generated
by $b_1^+,a^+,k'_1,k'_2,y'$ in  $\Frac\mathcal U(\mathfrak b^+)$ is isomorphic to the algebra $\mathbf A_1\otimes\mathbf S_3$.

We have now to formulate the commutation relations of the last generator $t'$ with the generators of $\mathcal V$.
It is clear that $t'b_1^+=b_1^+t'$ and $t'k'_1=k'_1t'$; we compute
$t'a^+-a^+t'=1$ and $t'k'_2=k'_2t'-(a^+)^{-1}k'_2+(a^+)^{-1}t'$. We try to replace $t'$ by
a generator of the form $a^+t'+p$ commuting with $a^+$ and $k'_2$, with $p\in\mathcal W$.
A solution is given by $p=-a^+k'_2$. In other words, the element 
$u:=a^+t'-a^+k'_2$ satisfies $[u,b_1^+]=[u,k'_1]=[u,a^+]=[u,k'_2]=0$. We calculate:\begin{align*}
uy' &= a^+t'y'-a^+k'_2y'=a^+t'(a^+)^{-1}y-y'a^+k'_2\\
    &=a^+((a^+)^{-1}t'-(a^+)^{-2})y-y'a^+k'_2=t'y-(a^+)^{-1}y-y'a^+k'_2\\
    &=yt'-y'-y'a^+k'_2=a^+y't'-y'a^+k'_2-y'\\
    &=y'(a^+t'-a^+k'_2)-y'=y'u-y'.\end{align*}
This relation becomes $y't''-t''y'=1$
with notation:\begin{equation}\label{deft''}t'':=(y')^{-1}u=(y')^{-1}a^+(b_1^+)^{-2}{t} - (y')^{-1}a^+k'_2.\end{equation}
Since $u$ commutes with $a^+, k'_2, b_1^+, k'_1$ it follows from \eqref{bracy'} that
$t''$  commutes with $a^+$ and $k'_2$, and anticommutes with $b_1^+$ and $k'_1$.

To sum up, starting from the generators $b_1^+, a^+, k_2,k_1, b_2^+, t$ of 
$\mathcal U(\mathfrak b^+)$, we have proved that the elements $b_1^+, a^+, k'_2, k'_1, y',t''$
defined by \eqref{defk'}, \eqref{defy'}, \eqref{deft''}  generate $\Frac\mathcal U(\mathfrak b^+)$.
The subalgebra generated by $k'_2$ and $a^+$ is isomorphic to the Weyl algebra $\mathbf A_1$,
the subalgebra  generated by $k'_1,b_1^+,y'$ and $t''$ is isomorphic to the algebra $\mathbf S_4$,
each element of the first subalgebra commutes with each element of the second one, and 
$\Frac\mathcal U(\mathfrak b^+)=\Frac(\mathbf A_1\otimes\mathbf S_4)$. Hence the proof is complete.\end{proof}

\subsection{Remark} The enveloping algebra $\mathcal U(\mathfrak b_{\overline 0}^+)$ is the subalgebra of $\mathcal U(\mathfrak b^+)$ 
generated by $(b_1^+)^2,(b_2^+)^2,a^+, t,k_1,k_2$. With the notations used in the proof of theorem \ref{thmB}, 
the generators $a^+,k'_2$ lie in $\mathcal U(\mathfrak b_{\overline 0}^+)$ and we define  in $\mathcal U(\mathfrak b_{\overline 0}^+)$
the elements $\ell_1:=(b_1^+)^{-2}({k_1}-k_2)$, $y'':=(y')^2=-(a^+)^{-2}(b_1^+)^2({b_2^+})^2+1$ and
$t''':=\frac12(y'')^{-1}a^+(b_1^+)^{-2}{t} - \frac12(y'')^{-1}a^+k'_2$. Then 
$\Frac\mathcal U(\mathfrak b_{\overline 0}^+)$ is generated by $k'_2,a^+,(b_1^+)^2,\ell_1,y'',t'''$ and
the brackets between these generators are $[k'_2,a^+]=[\ell_1,(b_1^+)^2]=[y'',t''']=1$ 
and 0 in all other cases.
We recover the well known Gelfand-Kirillov property that $\Frac\mathcal U(\mathfrak b_{\overline 0}^+)$
is a classical Weyl skew field $D_3$ over a trivial center $\C$.\bigskip

The following theorem gives a decomposition of $\Frac\mathcal U(\mathfrak p^+)$
into two commuting subfields respectively isomorphic to  $\mathbf D_2^0$ and $\mathbf F_3$.\medskip

\subsection{Theorem}\label{thmP} {\sl $\Frac\mathcal U(\mathfrak p^+)$ is isomorphic to 
$\Frac(\mathbf A_1\otimes\mathbf A_1\otimes\mathbf S_3)$.}

\begin{proof} By \ref{P}, $\mathcal U(\mathfrak p^+)$ is generated in $\mathcal U(\mathfrak{osp}(1,4))$
by $\mathcal U(\mathfrak b^+)$ and $b_2^-$ with commutation relations coming from \eqref{bracB+}, \eqref{bracP+}
and \eqref{eqVn+}. We start replacing in $\Frac\mathcal U(\mathfrak n^+)$ the generators $t$ and $a^+$ 
by:\begin{equation}\label{defuv}u_1:=t,\ \ v_1:=a^+(b_1^+)^{-2},\end{equation}
which commute with $b_1^+$ and satisfy $u_1v_1-v_1u_1=1$. Then we consider the enveloping algebra $\mathcal U(\mathfrak l)$
of the Levi subalgebra generated by $b_2^+,b_2^-,k_2$, see \ref{L}. In $\Frac\mathcal U(\mathfrak p^+)$, we replace
$b_2^+$ and $b_2^-$ by $m_2^+:=b_2^+-a^+(b_ 1^+)^{-1}$ {and} $m_2^-:=b_2^--t(b_ 1^+)^{-1}$ in order to simplify the commutation
relations with the previous generators $b_1^+,u_1,v_1$:\begin{equation}\label{m2}
[m_2^{\pm},u_1]=[m_2^{\pm},v_1]=0\qquad\text{et}\qquad m_2^\pm b_1^+=-b_1^+ m_2^\pm.
\end{equation} We define $\ell_2:=\frac12(m_2^+m_2^-+m_2^-m_2^+)$. A technical but straightforward calculation 
gives $\ell_2=k_2-ta^+(b_1^+)^{-2}+{\frac12}=k_2-u_1v_1+\frac12$.
Since $m_2^+m_2^-$ and $m_2^-m_2^+$ commute with $b_1^+,u_1,v_1$ by \eqref{m2}, the same is true for $\ell_2$. Moreover we compute:
$[\ell_2,m_2^+]=m_2^+$ and $[\ell_2,m_2^-]=-m_2^-$. 
The subalgebra generated by $m_2^+,m_2^-,\ell_2$ is isomorphic
to $\mathcal U(\mathfrak{osp}(1,2))$ and we apply the method used
in proposition \ref{mainresultosp12} setting:\begin{equation}\label{defuv2}
u_2:=-(m_2^-)^{-1}\ell_2, \ \  v_2:= m_2^-, \ \ z_2:=-2m_2^-m_2^++2\ell_2+1.\end{equation}
To sum up, the subfield $L$ of $\Frac\mathcal U(\mathfrak p^+)$ 
generated by $b_1^+, t, a^+,b_2^-,k_2,b_2^+$ is also generated by $b_1^+,u_1,v_1,v_2,z_2,u_2$
with relations:{\renewcommand{\arraystretch}{1.5}\begin{equation*}\begin{matrix}
[u_1,v_1]=1, \hfill&[u_2,v_2]=1,\hfill   &[u_1,v_2]=[u_2,v_1]=[u_1,u_2]=[v_1,v_2]=0,\hfill \\
b_1^+u_2=-u_2b_1^+,\hfill   &b_1^+v_2=-v_2b_1^+,\hfill  &[b_1^+,u_1]=[b_1^+,v_1]=[b_1^+,z_2]=0, \hfill \\
z_2u_2=-u_2z_2,\hfill  &z_2v_2=-v_2z_2,\hfill   &[z_2,u_1]=[z_2,v_1]=0. \hfill\\ 
\end{matrix}\end{equation*}}We can replace the generator $b_1^+$ by $w_1:=z_2^{-1}b_1^+$ which is central in $L$.\medskip

In the last step we look at the action of $k_1$ on $L$. Technical calculations using \eqref{defuv} and \eqref{defuv2}
show that, on one hand $[k_1,u_2]=[k_1,v_2]=[k_1,z_2]=0$, and on the other hand
$[k_1,u_1]=u_1,[k_1,v_1]=-v_1, [k_1,w_1]=w_1$. As in lemma 4 of \cite{AOVdB}, the last
change of variable $k''_1:=(k_1+u_1v_1)w_1^{-1}$ doesn't change the first three relations
and changes the last three into: $[k''_1,u_1]=[k''_1,v_1]=0$ and $[k''_1,w_1]=1$.
We conclude that in $\Frac\mathcal U(\mathfrak p^+)$ 
the subalgebra generated by $k''_1,w_1$ is isomorphic to $\mathbf A_1$,
the subalgebra generated by $u_1,v_1$ is  also isomorphic to $\mathbf A_1$,
the subalgebra generated by $u_2,v_2,z_2$ is isomorphic to $\mathbf S_3$, 
and $\Frac\mathcal U(\mathfrak p^+)$ is isomorphic to 
$\Frac(\mathbf A_1\otimes\mathbf A_1\otimes\mathbf S_3)$.\end{proof}

\subsection{Remark}\label{L79}The enveloping algebra $\mathcal U(\mathfrak p_{\overline 0}^+)$ is the subalgebra of $\mathcal U(\mathfrak p^+)$ 
generated by $(b_1^+)^2,(b_2^+)^2,a^+, t,k_1,k_2,(b_2^-)^2$. 
Computing the brackets between these generators we find exactly the table of the Lie algebra
denoted by $L_{7,9}$ in \cite{AOVdB} p.\,565 up to the following change of variables:
{\renewcommand{\arraystretch}{1.5}\begin{equation*}\begin{matrix}
&\textstyle e_{0}:=t, \hfill &e_1:=a^+, \hfill  &e_2:=(b_1^+)^2, \hfill  &e_3:=k_1, \\
&x:=\frac12(b_2^-)^2,\hfill &y:=-\frac12(b_2^+)^2,  \hfill &h:=-k_2.\hfill &
\end{matrix}\end{equation*}}It is proved in \cite{AOVdB} that the Lie algebra 
$\mathfrak p_{\overline 0}^+=L_{7,9}$
satisfies the Gelfand-Kirillov property with
$\Frac\mathcal U(\mathfrak p_{\overline 0}^+)=\Frac(\mathbf A_1\otimes\mathbf A_1\otimes\mathbf A_1\otimes\C[c])$.
The central generator $c$ and the pairs of elements $p_i,q_i$ ($1=1,2,3$) described in \cite{AOVdB}
as generators of each copy of $\mathbf A_1$ 
correspond with our notations in the proof of theorem \ref{thmP} to:
{\renewcommand{\arraystretch}{1.5}\begin{equation*}\label{lexic}\begin{matrix}
\textstyle 
p_1:=u_1,\hfill & \ \ p_2:=v_ 2^{-1}u_2,\hfill & \ \ p_3:=\frac12 k''_1w_1^{-1}z_2^{-2},\hfill & \ \ c:=\frac14(z_2+1)^2-1,\\
q_1:=v_1,\hfill & \ \ q_2:=\frac12(v_2)^2,\hfill & \ \ q_3:=w_1^2z_2^2,\hfill &\\
\end{matrix}\end{equation*}}which gives an explicite description of the embedding:
\[\Frac\mathcal U(\mathfrak p_{\overline 0}^+)=\Frac(\mathbf A_1\otimes\mathbf A_1\otimes\mathbf A_1\otimes\C[c])\subset
\Frac\mathcal U(\mathfrak p^+)=\Frac(\mathbf A_1\otimes\mathbf A_1\otimes\mathbf S_3).\] 

\subsection{Illustration}With the conventions of remark \ref{graph}, proposition \ref{propN}, theorem
\ref{thmB} and theorem \ref{thmP} can be represented by the following pictures:\bigskip

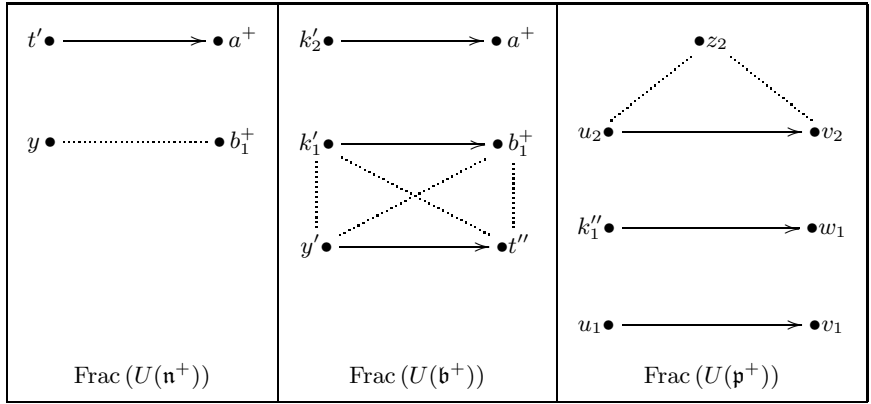
\centerline{\footnotesize{\renewcommand{\arraystretch}{2.1}\begin{tabular}{|c|c|c|}
\hline
$\xymatrix{
t' \bullet  \ar@*{[|<2pt>]}[rr]& & {\bullet}\,a^+ \\
y\, \bullet  \ar@{..}[rr]& & {\bullet}\,b_1^+}$&
$\xymatrix{
k'_2 \bullet \ar[rr] & & {\bullet}\,a^+ \\
k'_1 \bullet \ar[rr] & & {\bullet}\,b_1^+\\
y' \bullet\ar@{..}[rru] \ar[rr] \ar@{..}[u] & &\ar@{..}[llu]  \ar@{..}[u]  \bullet t''\\ }$&
$\xymatrix{
  &   \ar@{..}[ld] \bullet z_2\ar@{..}[rd]& \\
u_2 \bullet  \ar@*{[|<2pt>]}[rr]& & {\bullet}v_2 \\
k''_1 \bullet  \ar@*{[|<2pt>]}[rr]& & {\bullet}w_1   \\
u_1 \bullet  \ar@*{[|<2pt>]}[rr]& & {\bullet}v_1}$ \\
$\Frac(U(\mathfrak n^+))$ &$\Frac(U(\mathfrak b^+))$ &$\Frac(U(\mathfrak p^+))$\\
\hline
\end{tabular}}}

\subsection{Corollary}{\sl  The center of $\Frac\mathcal U(\mathfrak n^+)$ is 
a purely transcendental extension of $\C$ of degree two,
the center of $\Frac\mathcal U(\mathfrak b^+)$ is $\C$, and 
the center of $\Frac\mathcal U(\mathfrak p^+)$ is a purely transcendental extension of $\C$ of degree one.}
\begin{proof} Follows directly from proposition \ref{propN}, theorem \ref{thmB} and theorem \ref{thmP}
applying the results  on the centers \ref{propD}.(ii) and \ref{propF3F4}.(ii).
More explicitly with the notations used in the proofs, the center of $\Frac\mathcal U(\mathfrak n^+)$
is $\C((b_1^+)^2 ,y^2 )$ and the center of $\Frac\mathcal U(\mathfrak p^+)$ is $\C(z_2^2)$.\end{proof}

\subsection{Proposition}{\sl The skew fields $\Frac\mathcal U(\mathfrak b^+)$ and $\Frac\mathcal U(\mathfrak p^+)$
are not isomorphic to $\mathbf D_{r,t}^s$ for any $r,s,t\geq 0$.}
\begin{proof}Suppose that $\Frac\mathcal U(\mathfrak b^+)$ is isomorphic to some
skew field $\mathbf D_{r,t}^s$. Comparing the Gelfand-Kirillov transcendence degrees and the centers,
we have $2r+2s+t=6$ and $2s+t=0$, hence $\Frac\mathcal U(\mathfrak b^+)$ would be isomorphic 
to the usual Weyl skew field $\mathbf D_{3,0}^0=D_3(\C)$ which is impossible
because, as at the end of the proof of proposition \ref{propF3F4}, we have 
$G(D_3(\C))=\{1\}$ and $-1\in G(\Frac\mathcal U(\mathfrak b^+))$. 
Suppose now that $\Frac\mathcal U(\mathfrak p^+)$ is isomorphic to some
skew field $\mathbf D_{r,t}^s$. We obtain $2r+2s+t=7$ 
and $2s+t=1$, hence $\Frac\mathcal U(\mathfrak p^+)$ would be isomorphic 
to $\mathbf D_{3,1}^0$, which is impossible by the same argument.\end{proof}

\subsection{Remark}\label{L77} The Lie algebra ${\mathfrak{sp}}(4)$ contains two non isomorphic parabolic subalgebras 
corresponding to the cases denoted by $L_{7,7}$ and $L_{7,9}$ in the classification of \cite{AOVdB}.
We have seen in \ref{L79} that the even part $\mathfrak p_{\overline 0}^+$ of the parabolic subsuperalgebra $\mathfrak p^+$
is isomorphic to $L_{7,9}$. But we can also define a subsuperalgebra $\mathfrak q^+$ of $\mathfrak{osp}(1,4)$ whose even part is the alternative
parabolic subalgebra $L_{7,7}$ of ${\mathfrak{sp}}(4)$. It is defined by $\mathfrak q^+=\mathfrak q_{\overline 0}^+\oplus\mathfrak g_{\overline 1}^+$ with $\mathfrak g_{\overline 1}^+=\C b_1^+\oplus\C b_2^+$ and $\mathfrak q_{\overline 0}^+=\mathfrak b_{\overline 0}^+\oplus \C s$,
where $s$ is defined in \eqref{gen}. A basis of $\mathfrak q_{\overline 0}^+$ is $\{c_1^+, c_2^+, a^+, t, k_1,k_2,s\}$ and
computing the brackets in $\mathfrak{osp}(1,4)$ we retrieve the table of $L_{7,7}$ in \cite{AOVdB} up to the following change
of notations:   
{\renewcommand{\arraystretch}{1.5}\begin{equation*}\begin{matrix}
&\textstyle e_{0}:=c_1^+, \hfill &e_1:=2a^+, \hfill  &e_2:=c_2^+, \hfill  &e_3:=-\frac12(k_1+k_2), \\
&x:=t,\hfill &y:=s,  \hfill &h:=k_1-k_2.\hfill &
\end{matrix}\end{equation*}}By a method similar to that of theorem \ref{thmP}, we can prove that $\Frac\mathcal U(\mathfrak q^+)$ 
is also isomorphic to 
$\Frac(\mathbf A_1\otimes\mathbf A_1\otimes\mathbf S_3)$.

\section*{Acknowledgements} We would like to thank Alfons Ooms for drawing our attention to the case 
of the second parabolic subalgebra considered in remark \ref{L77}.

\end{document}